\documentclass[11pt,a4paper]{article}

\usepackage{inputenc}
\usepackage{amsmath}
\usepackage{bm}
\usepackage{bbold}
\usepackage{amsthm}

\usepackage{curves}

\usepackage{hyperref}

\title{Evaluation of the Mean Cycle Time\\ in Stochastic Discrete Event Dynamic Systems\thanks{Proc. 6th Intern. Conf. on Queueing Theory and Network Applications (QTNA'11), ACM, New York, 2011, pp.~93--100.}
} 

\author{Nikolai Krivulin\thanks{Faculty of Mathematics and Mechanics, St.~Petersburg State University, 28 Universitetsky Ave., St.~Petersburg, 198504, Russia, 
nkk@math.spbu.ru.} \thanks{The work was partially supported by the Russian Foundation for Basic Research under Grant \#09-01-00808.}
}

\date{}

\begin{document}

\maketitle

\begin{abstract}
We consider stochastic discrete event dynamic systems that have time evolution represented with two-dimensional state vectors through a vector equation that is linear in terms of an idempotent semiring. The state transitions are governed by second-order random matrices that are assumed to be independent and identically distributed. The problem of interest is to evaluate the mean growth rate of state vector, which is also referred to as the mean cycle time of the system, under various assumptions on the matrix entries. 

We give an overview of early results including a solution for systems determined by matrices with independent entries having a common exponential distribution. It is shown how to extend the result to the cases when the entries have different exponential distributions and when some of the entries are replaced by zero. Finally, the mean cycle time is calculated for systems with matrices that have one random entry, whereas the other entries in the matrices can be arbitrary nonnegative and zero constants. The random entry is always assumed to have exponential distribution except for one case of a matrix with zero row when the particular form of the matrix makes it possible to obtain a solution that does not rely on exponential distribution assumptions.
\\

\end{abstract}

\section{Introduction}

Many actual systems in engineering, manufacturing, management, and other areas can be modeled as stochastic discrete event dynamic systems. The dynamics of the systems is normally described by recurrent equations for a set of state variables representing the occurrence time of system events. In some cases, the state evolution of the system can also be written in vector form with a state transition matrix through a vector equation that is linear in the sense of an idempotent semiring \cite{Kolokoltsov1997Idempotent,Golan2003Semirings,Heidergott2006Maxplus,Heidergott2007Maxplus}. As an example, one can consider the G/G/1 queue and related representation in the form of Lindley recursions (see, e.g. \cite{Krivulin2009Evaluation}).

In the analysis of system performance, one is often interested in evaluating the mean growth rate of the system state vector, which presents the mean time of one system operation cycle. However, exact calculation of the mean cycle time can be a rather cumbersome problem even for quite simple systems. Most of the known results are limited to systems with second-order state transition matrices having independent random entries.

In this paper we give an overview of existing solutions and present new results for the problem. We start with motivating examples of calculating the mean cycle time drawn from manufacturing and telecommunications. Early results from \cite{Resing1990Asymptotic,Olsder1990Discrete,Jean-marie1994Analytical} including solutions for second-order matrices with independent and identically distributed (i.i.d.) entries are outlined. For exponential, continuous uniform, Bernoulli, geometric, or discrete uniform distributions, the solutions are given as rational functions of distribution parameters.

To get solutions for other systems with second-order matrices, we implement an approach based on constructing and analyzing a sequence of one-dimensional distributions. The mean cycle time is calculated as the mean of a random variable (r.v.) determined by the limiting distribution of the sequence. Specifically, the case of a matrix where the entries have exponential distributions with arbitrary parameters is considered. We describe a computational technique \cite{Krivulin2008Evaluation} that actually reduces the problem to the solution of an appropriate set of linear algebraic equations. 

Furthermore, systems with random matrices that can have some constant entries are examined. Both recent results for matrices with zero entries \cite{Krivulin2007Growth,Krivulin2009Calculating} and new results for matrices that can have arbitrary nonnegative constant entries are presented. The new results include solutions for systems with second-order matrices that have one random entry, whereas the other entries can be arbitrary nonnegative and zero constants. The random entry is always assumed to have exponential distribution except for one case when the particular form of the matrix makes it possible to obtain a solution that does not rely on exponential distribution assumptions.

\section{Preliminary Observations}

\subsection{Motivating Examples}

We start with example systems drawn from manufacturing and telecommunications. Other related examples can be found in \cite{Resing1990Asymptotic,Jean-marie1994Analytical,Heidergott2006Maxplus}.

\subsubsection{Manufacturing}

Consider a manufacturing system that consists of two production centers $A$ and $B$. Each center produces its output based on use of an output from the other center. The operation of each center as well as of the entire system consists in a sequence of production cycles. Every cycle involves simultaneous production of a new output and transportation of the output by the previous cycle. A center completes its current cycle as soon as it finishes production of the current output and an output from the other center arrives. A current production cycle of the system comes to the end and the next cycle is initiated as soon as both centers $A$ and $B$ complete their respective production cycles.

Suppose that at the initial time, each center has its output available for delivery. Given random production and transportation time for both centers, one is often interested in evaluating the mean time of production cycle as the number of cycles goes to infinity. The system throughput defined as the inverse of the mean cycle time can be another performance measure of interest.

For every cycle $k=1,2,\ldots$ we introduce the following notation. Let $\alpha_{k}$ and $\delta_{k}$ be r.v.'s that represent the production time at respective centers $A$ and $B$, and let $\beta_{k}$ and $\gamma_{k}$ be r.v.'s that represent the transportation time from $B$ to $A$ and from $A$ to $B$. By $x(k)$ and $y(k)$ we denote the cycle completion epochs at centers $A$ and $B$. 

With the condition that $x(0)=y(0)=0$, the dynamics of the system is naturally represented by two equations
\begin{align}
\begin{split}
x(k)
&=
\max(x(k-1)+\alpha_{k},y(k-1)+\beta_{k}),
\\
y(k)
&=
\max(x(k-1)+\gamma_{k},y(k-1)+\delta_{k}).
\end{split}
\label{E-xkyk}
\end{align}

The mean cycle time for the system is given by
\begin{equation}
\lambda
=
\lim_{k\to\infty}\frac{1}{k}\max(x(k),y(k)),
\label{E-lambda}
\end{equation}
provided that the limit on the right exists.

\subsubsection{Telecommunications}

Suppose that there is a system of two work stations $A$ and $B$ which exchange messages over a communication network. Each station generates and sends messages in response to messages received from the other station. If a message is sent from a station, it immediately enters the network and then goes from one intermediate node to another until it arrives at the other station.

The operation of each station and of the entire system consists in a sequence of communication sessions. A new session at a station begins with generation of a new message and completes by sending the message to the other station. Once the message is generated, it is sent from the station as soon as a message from the other station arrives. A current communication session of the system starts when both nodes $A$ and $B$ complete their previous sessions, and lasts until they complete the current sessions.

For every session $k=1,2,\ldots$ we denote by $x(k)$ and $y(k)$ the session completion epochs for respective stations $A$ and $B$. Let $\alpha_{k}$ and $\delta_{k}$ be the random message generation time for stations $A$ and $B$, and $\beta_{k}$ and $\gamma_{k}$ be the random message transmission time from $B$ to $A$ and from $A$ to $B$.

Given the r.v.'s $\alpha_{k}$, $\beta_{k}$, $\gamma_{k}$, and $\delta_{k}$ for all $k=1,2,\ldots$ we consider the problem of evaluating the mean communication session time which one can take as a natural performance measure for the system. As it is easy to see, the dynamic equations for the system have the same form as \eqref{E-xkyk}. Moreover, the mean communication session time can be given by \eqref{E-lambda} as well.

\subsection{Problem Representation and Preliminary Analysis}

Consider an idempotent semiring with two operations, maximum as semiring addition, and arithmetic addition as semiring multiplication. Define a system state transition matrix for every $k=1,2,\ldots$ as
\begin{equation}
A(k)
=
\left(
	\begin{array}{cc}
		\alpha_{k}	& \beta_{k} \\
		\gamma_{k}	& \delta_{k}
	\end{array}
\right).
\label{E-Ak}
\end{equation}

We suppose that each of the sequences $\{\alpha_{k}\}$, $\{\beta_{k}\}$, $\{\gamma_{k}\}$, and $\{\delta_{k}\}$ involves either i.i.d. r.v.'s or nonnegative constants, and that the sequences are independent themselves.

Furthermore, we introduce a state vector $\bm{z}(k)$ as follows
$$
\bm{z}(k)
=
\left(
	\begin{array}{c}
		x(k) \\
		y(k)
	\end{array}
\right),
\qquad
\bm{z}(0)
=
\left(
	\begin{array}{c}
		0 \\
		0
	\end{array}
\right).
$$

Scalar equations \eqref{E-xkyk} can now be replaced with an equation
\begin{equation}
\bm{z}(k)
=
A(k)\bm{z}(k-1),
\label{E-zk}
\end{equation}
where matrix-vector multiplication is thought of in the sense of addition and multiplication in the semiring.

With the norm $\|\cdot\|$ defined as the maximum entry of a vector, we can rewrite the limit in \eqref{E-lambda} to represent the problem of interest as that of finding
\begin{equation}
\lambda
=
\lim_{k\to\infty}\|\bm{z}(k)\|^{1/k},
\label{E-lambda1}
\end{equation}
where $\lambda$ can generally be interpreted as the mean growth rate of state vector in system \eqref{E-zk}.

Consider a stochastic system governed by equation \eqref{E-zk} and assume that the random entries in matrix \eqref{E-Ak} have finite means. It is not difficult to verify (see e.g. \cite{Kingman1973Subadditive,Krivulin2009Evaluation}) that under this assumption, the limit at \eqref{E-lambda1} exists with probability one. Moreover, in this case, it holds that
$$
\lambda
=
\lim_{k\to\infty}\mathsf{E}\|\bm{z}(k)\|^{1/k},
$$
or, in ordinary notation, that
$$
\lambda
=
\lim_{k\to\infty}\frac{1}{k}\mathsf{E}\max(x(k),y(k)).
$$ 

Note that the above vector representation, among other benefits, enables one to enhance the analysis of the system by using models and methods of idempotent algebra (see e.g., \cite{Krivulin2009Evaluation}). Specifically, the range of particular cases to analyze can be reduced on the basis of the following observations. First, it is not difficult to see from \eqref{E-zk} and \eqref{E-lambda1} that both systems, with a matrix $A(k)$ and its transpose $A^{T}(k)$, have the same value of $\lambda$. Moreover, the mean cycle time is invariant to simultaneous row and column permutations in $A(k)$. Finally, any superposition of the above transformations also leaves $\lambda$ unchanged. 

\subsection{Previous Results}

The problem of evaluating the mean cycle time may involve cumbersome computations even for quite simple systems. Many known results are obtained for systems with second-order matrices with i.i.d. random entries \cite{Olsder1990Discrete,Resing1990Asymptotic,Jean-marie1994Analytical}. 

A system with a matrix where all entries are exponentially distributed with parameter $\mu$ is examined in \cite{Olsder1990Discrete,Resing1990Asymptotic}. A somewhat too complicated solution based on the theory of Markov chains in \cite{Resing1990Asymptotic} leads to the result
$\lambda=407/(228\mu)$.

For the particular case $\mu=1$, this result is obtained in \cite{Olsder1990Discrete} with a less rigorous approach that actually reduces the problem to the solution of an integral equation.

Furthermore, the technique developed in \cite{Resing1990Asymptotic} is extended to get the solution $\lambda=0.719$ when all entries of the matrix have a continuous uniform distribution on $[0,1]$. 

Another rather sophisticated approach to the problem with second-order matrices is proposed in \cite{Jean-marie1994Analytical}. The approach considers a sequence of two-dimensional distributions of the system state vector and examines the convergence of the sequence by applying the Laplace transform of the distributions. It is shown that for a system with a matrix where the diagonal and off-diagonal entries are exponentially distributed with respective parameters $\mu$ and $\nu$, the solution takes the form
$$
\lambda=P(\mu,\nu)/Q(\mu,\nu),
$$
where
\begin{align*}
P(\mu,\nu)
&=
160\mu^{10}+1776\mu^{9}\nu+8220\mu^{8}\nu^{2}+21378\mu^{7}\nu^{3}
\\
&\phantom{=}
+35595\mu^{6}\nu^{4}+41566\mu^{5}\nu^{5}+35595\mu^{4}\nu^{6}
\\
&\phantom{=}
+21378\mu^{3}\nu^{7}+8220\mu^{2}\nu^{8}+1776\mu\nu^{9}+160\nu^{10},
\\
Q(\mu,\nu)
&=
16\mu\nu(\mu+\nu)(8\mu^{8}+80\mu^{7}\nu+321\mu^{6}\nu^{2}+690\mu^{5}\nu^{3}
\\
&\phantom{=}
+880\mu^{4}\nu^{4}+690\mu^{3}\nu^{5}+321\mu^{2}\nu^{6}+80\mu\nu^{7}+8\nu^{8}).
\end{align*}

For a system with a symmetric matrix having both diagonal entries exponentially distributed with parameter $\mu$ and off-diagonal entries equal to zero, the solution is $\lambda=5/(4\mu)$.

We conclude with an overview of results for discrete distributions. In \cite{Olsder1990Discrete,Resing1990Asymptotic}, the case of matrices of order $n=2,3$ is examined when the entries have the discrete uniform distribution on $0,1,\ldots,m$. The evolution of the differences between components of the state vector is represented as a Markov chain, whereas the mean cycle time is evaluated from the stationary probabilities of the chain. The results obtained can be summarized as follows
\begin{align*}
&n=2,\quad m=1,	& \lambda&=6/7\approx0.857,
\\
&n=2,\quad m=2,	& \lambda&=0.803,
\\
&n=2,\quad m\to\infty,	& \lambda&\to0.719,
\\
&n=3,\quad m=1,	& \lambda&=0.979.
\end{align*}

The above approach is further developed in \cite{Jean-marie1994Analytical} where a technique is implemented based on recursive computation of joint generating function of the state vector. In the case of a matrix with i.i.d. entries having Bernoulli distribution with parameter $p$, the solution is given by   
$$
\lambda
=
1
-
\frac{(1+2p)(1-p)^{4}}{1-2p(1-p)(1-3p+p^{2})}.
$$

Finally, for the system where the entries of the matrix are geometrically distributed with parameter $p$, it holds that
$$
\lambda
=
N(p)/D(p),
$$
where
{\thinmuskip=1mu\medmuskip=1mu plus 1mu minus 1mu\thickmuskip=3mu plus 3mu%
\begin{align*}
N(p)
&=
p(4+18p+50p^{2}+99p^{3}+175p^{4}+244p^{5}+289p^{6}
\\
&\phantom{=}
+273p^{7}+218p^{8}+137p^{9}+77p^{10}+32p^{11}+11p^{12}+p^{13}),
\\
D(p)
&=
(1-p)(1+p)(1+p+p^{2})(1+6p+8p^{2}+20p^{3}
\\
&\phantom{=}
+25p^{4}+32p^{5}+25p^{6}+20p^{7}+8p^{8}+6p^{9}+p^{10}).
\end{align*}
}

\section{An Overview of Recent Results}

This section offers solutions that generalize results given in \cite{Olsder1990Discrete,Resing1990Asymptotic,Jean-marie1994Analytical} for systems that have second-order matrices with exponentially distributed entries. We start with a solution developed in \cite{Krivulin2008Evaluation} for the case of pure random matrices having only random entries with different exponential distributions. 

Furthermore, we consider results obtained in \cite{Krivulin2007Growth,Krivulin2009Calculating} for matrices with some nonrandom entries that are equal to zero. The reduced number of random entries in the matrices allows one to simplify the evaluation of the mean cycle time and usually gives quite compact results.

\subsection{The System With a Pure Random Matrix}

Consider a system with the state transition matrix
$$
A(k)
=
\left(
\begin{array}{cc}
	\alpha_{k}	& \beta_{k} \\
	\gamma_{k}	& \delta_{k}
\end{array}
\right),
$$
where $\alpha_{k}$, $\beta_{k}$, $\gamma_{k}$, and $\delta_{k}$ are assumed to have exponential distributions with respective parameters $\mu$, $\nu$, $\sigma$, and $\tau$.

To get the mean cycle time, a computational technique proposed in \cite{Krivulin2008Evaluation} reduces the problem to the solution of a system of linear equations, accompanied by the evaluation of a linear functional of the solution. The technique leans upon construction and examination of a sequence of probability density functions. It is shown that there is a one-to-one correspondence between the sequence of density functions and a sequence of certain vectors in a vector space. The correspondence is then exploited to provide for the result in terms of algebraic computations. 

Based on the above technique, the mean cycle time is evaluated as follows. First we introduce the vectors
$$
\bm{\omega}_{1}
=
(\omega_{10},\omega_{11},\omega_{12},\omega_{13})^{T},
\qquad
\bm{\omega}_{2}
=
(\omega_{20},\omega_{21},\omega_{22},\omega_{23})^{T}.
$$

Then we put
$$
\bm{\omega}
=
\left(
\begin{array}{c}
\bm{\omega}_{1} \\
\bm{\omega}_{2}
\end{array}
\right),
\qquad
W
=
\left(
\begin{array}{cc}
U_{1}V_{11} & U_{1}V_{12} \\
U_{2}V_{21} & U_{2}V_{22}
\end{array}
\right),
$$
where the matrices $U_{1}$ and $U_{2}$ are defined as
\begin{align*}
U_{1}
&=
\left(
\begin{array}{ccc}
1 & 1 & 1 \\
\frac{\mu}{\mu+\nu} & \frac{1}{2} & \frac{\mu+\nu}{\mu+2\nu} \\
\frac{\mu}{\mu+\tau} & \frac{\nu}{\nu+\tau} & \frac{\mu+\nu}{\mu+\nu+\tau} \\
\frac{\mu}{\mu+\nu+\tau} & \frac{\nu}{2\nu+\tau} & \frac{\mu+\nu}{\mu+2\nu+\tau}
\end{array}
\right),
\\
U_{2}
&=
\left(
\begin{array}{ccc}
1 & 1 & 1 \\
\frac{\sigma}{\mu+\sigma} & \frac{\tau}{\mu+\tau} & \frac{\sigma+\tau}{\mu+\sigma+\tau} \\
\frac{1}{2} & \frac{\tau}{\sigma+\tau} & \frac{\sigma+\tau}{2\sigma+\tau} \\
\frac{\sigma}{\mu+2\sigma} & \frac{\tau}{\mu+\sigma+\tau} & \frac{\sigma+\tau}{\mu+2\sigma+\tau}
\end{array}
\right).
\end{align*}

Furthermore, the matrices $V_{11}$ and $V_{12}$ are defined as
\begin{align*}
V_{11}
&=
\left(
\begin{array}{cccc}
\frac{\sigma}{\mu+\sigma} & 0 & -\frac{\mu\sigma}{(\mu+\tau)(\mu+\sigma+\tau)} & 0 \\
0 & \frac{\sigma}{\nu+\sigma} & 0 & -\frac{\nu\sigma}{(\nu+\tau)(\nu+\sigma+\tau)} \\
0 & -\frac{\sigma}{\mu+\nu+\sigma} & 0 & \frac{\sigma(\mu+\nu)}{(\mu+\nu+\tau)(\mu+\nu+\sigma+\tau)}
\end{array}
\right),
\\
V_{12}
&=
\left(
\begin{array}{cccc}
0 & \frac{\tau}{\mu+\tau} & 0 & -\frac{\mu\tau}{(\mu+\sigma)(\mu+\sigma+\tau)} \\
\frac{\tau}{\nu+\tau} & 0 & -\frac{\nu\tau}{(\nu+\sigma)(\nu+\sigma+\tau)} & 0 \\
0 & -\frac{\tau}{\mu+\nu+\tau} & 0 & \frac{\tau(\mu+\nu)}{(\mu+\nu+\sigma)(\mu+\nu+\sigma+\tau)}
\end{array}
\right),
\end{align*}
whereas $V_{21}$ and $V_{22}$ take the form
\begin{align*}
V_{21}
&=
\left(
\begin{array}{cccc}
\frac{\mu}{\mu+\sigma} & -\frac{\mu\sigma}{(\nu+\sigma)(\mu+\nu+\sigma)} & 0 & 0 \\
0 & 0 & \frac{\mu}{\mu+\tau} & -\frac{\mu\tau}{(\nu+\tau)(\mu+\nu+\tau)} \\
0 & 0 & -\frac{\mu}{\mu+\sigma+\tau} & \frac{\mu(\sigma+\tau)}{(\nu+\sigma+\tau)(\mu+\nu+\sigma+\tau)}
\end{array}
\right),
\\
V_{22}
&=
\left(
\begin{array}{cccc}
0 & 0 & \frac{\nu}{\nu+\sigma} & -\frac{\nu\sigma}{(\mu+\sigma)(\mu+\nu+\sigma)} \\
\frac{\nu}{\nu+\tau} & -\frac{\nu\tau}{(\mu+\tau)(\mu+\nu+\tau)} & 0 & 0 \\
0 & 0 & -\frac{\nu}{\nu+\sigma+\tau} & \frac{\nu(\sigma+\tau)}{(\mu+\sigma+\tau)(\mu+\nu+\sigma+\tau)}
\end{array}
\right).
\end{align*}

Suppose the vector $\bm{\omega}$ is the solution of the system
\begin{align*}
(I-W)\bm{\omega}
&=
\bm{0}, \\
\omega_{10}+\omega_{20}
&=
1,
\end{align*}
where $I$ represents the identity matrix, $\bm{0}$ is zero vector.

The mean cycle time is then given by
$$
\lambda
=
\bm{q}_{1}^{T}\bm{\omega}_{1}+\bm{q}_{2}^{T}\bm{\omega}_{2},
$$
where $\bm{q}_{1}$ and $\bm{q}_{2}$ are vectors such that
\begin{align*}
\bm{q}_{1}
&=
\left(
\begin{array}{cc}
\frac{\mu^{2}+\mu\sigma+\sigma^{2}}{\mu\sigma(\mu+\sigma)} \\
\frac{\mu\sigma(\mu+2\nu+\sigma)}{\nu(\mu+\nu)(\nu+\sigma)(\mu+\nu+\sigma)} \\
\frac{\mu\sigma(\mu+\sigma+2\tau)}{\tau(\mu+\tau)(\sigma+\tau)(\mu+\sigma+\tau)} \\
-
\frac{\mu\sigma(\mu+2\nu+2\tau+\sigma)}{(\nu+\tau)(\mu+\nu+\tau)(\nu+\sigma+\tau)(\mu+\nu+\sigma+\tau)}
\end{array}
\right),
\\
\bm{q}_{2}
&=
\left(
\begin{array}{c}
\frac{\nu^{2}+\nu\tau+\tau^{2}}{\nu\tau(\nu+\tau)} \\
\frac{\nu\tau(2\mu+\nu+\tau)}{\mu(\mu+\nu)(\mu+\tau)(\mu+\nu+\tau)} \\
\frac{\nu\tau(\nu+2\sigma+\tau)}{\sigma(\nu+\sigma)(\sigma+\tau)(\nu+\sigma+\tau)} \\
-
\frac{\nu\tau(2\mu+\nu+2\sigma+\tau)}{(\mu+\sigma)(\mu+\nu+\sigma)(\mu+\sigma+\tau)(\mu+\nu+\sigma+\tau)}
\end{array}
\right).
\end{align*}

\subsection{Systems With Matrices With Zero Entries}

Suppose the time to perform certain operations in a system (e.g., production of output at a center or transmission of messages between stations) is sufficiently short as compared to that of the other operations. Under this assumption, one can usually set the duration of the operations equal to zero with negligible loss of accuracy, and then consider equation \eqref{E-zk} with a matrix that includes one or more zero entries. 

The solution method in this case is based on construction of a sequence of probability distribution functions. The convergence of the sequence is examined and the limiting distribution is derived as the solution of an integral equation. The mean cycle time is then evaluated as the expected value of a r.v. determined through the limiting distribution function. The results obtained are as follows.

\subsubsection{Matrices With Zero Off-Diagonal Entries}

The system state transition matrix together with the related result take the form 
$$
A(k)
=
\left(
	\begin{array}{cc}
		\alpha_{k}	& 0 \\
		0						& \delta_{k}
	\end{array}
\right),
\quad
\lambda
=
\frac{\mu^{4}+\mu^{3}\tau+\mu^{2}\tau^{2}+\mu\tau^{3}+\tau^{4}}{\mu\tau(\mu+\tau)(\mu^{2}+\tau^{2})}.
$$

With $\tau=\mu=1$ we have the mean cycle time $\lambda=1.25$ which coincides with that in \cite{Jean-marie1994Analytical}.

\subsubsection{Matrices With Zero Diagonal}

In this case, the matrix and the mean cycle time are represented as
$$
A(k)
=
\left(
	\begin{array}{cc}
		0						& \beta_{k} \\
		\gamma_{k}	& 0
	\end{array}
\right),
\quad
\lambda
=
\frac{4\nu^{2}+7\nu\sigma+4\sigma^{2}}{6\nu\sigma(\nu+\sigma)}.
$$

\subsubsection{Matrices With Zero Row or Column}

Provided that the second row in the matrix has only zero entries, we arrive at
$$
A(k)
=
\left(
	\begin{array}{cc}
		\alpha_{k}	& \beta_{k} \\
		0						& 0
	\end{array}
\right),
\qquad
\lambda
=
\frac{2\mu^{4}+7\mu^{3}\nu+10\mu^{2}\nu^{2}+11\mu\nu^{3}+4\nu^{4}}{\mu\nu(\mu+\nu)^{2}(3\mu+4\nu)}.
$$

Considering that $\lambda$ is invariant to transposition of the matrix, we immediately get the solution when the entries of the second column are zero
$$
A(k)
=
\left(
	\begin{array}{cc}
		\alpha_{k}	& 0 \\
		\gamma_{k}	& 0
	\end{array}
\right),
\qquad
\lambda
=
\frac{2\mu^{4}+7\mu^{3}\sigma+10\mu^{2}\sigma^{2}+11\mu\sigma^{3}+4\sigma^{4}}{\mu\sigma(\mu+\sigma)^{2}(3\mu+4\sigma)}.
$$

\subsubsection{Matrices With One Zero Diagonal Entry}

Consider a system with the state transition matrix
$$
A(k)
=
\left(
	\begin{array}{cc}
		\alpha_{k}	& \beta_{k} \\
		\gamma_{k}	& 0
	\end{array}
\right).
$$

Although evaluation of the mean cycle time for this system in the general case leads to rather cumbersome algebraic manipulations, there are two main particular cases that offer results in a relatively compact form. Under the condition that $\sigma=\mu$, we have
$$
\lambda
=
\frac{48\mu^{5}+238\mu^{4}\nu+495\mu^{3}\nu^{2}+581\mu^{2}\nu^{3}+326\mu\nu^{4}+68\nu^{5}}{2\mu\nu(36\mu^{4}+147\mu^{3}\nu+215\mu^{2}\nu^{2}+130\mu\nu^{3}+28\nu^{4})}.
$$

Provided that $\sigma=\nu$, the mean cycle time is given by
$$
\lambda=P(\mu,\nu)/Q(\mu,\nu),
$$
where
{\thinmuskip=2mu\medmuskip=3mu plus 2mu minus 3mu\thickmuskip=4mu plus 4mu%
\begin{align*}
P(\mu,\nu)
&=
15\mu^{8}+152\mu^{7}\nu+624\mu^{6}\nu^{2}+1382\mu^{5}\nu^{3}+1838\mu^{4}\nu^{4}
\\
&\phantom{=}
+1592\mu^{3}\nu^{5}+973\mu^{2}\nu^{6}+384\mu\nu^{7}+64\nu^{8},
\\
Q(\mu,\nu)
&=
\mu\nu(\mu+\nu)^{2}(12\mu^{5}+97\mu^{4}\nu+286\mu^{3}\nu^{2}+397\mu^{2}\nu^{3}
\\
&\phantom{=}
+256\mu\nu^{4}+64\nu^{5}).
\end{align*}
}

\subsubsection{Matrices With One Zero Off-Diagonal Entry}

Suppose that there is a system with a state transition matrix defined as
$$
A(k)
=
\left(
	\begin{array}{cc}
		\alpha_{k}	& \beta_{k} \\
		0						& \delta_{k}
	\end{array}
\right).
$$

Consider two particular cases. With the condition $\nu=\mu$, we have
$$
\lambda
=
P(\mu,\tau)/Q(\mu,\tau),
$$
where
{\thinmuskip=2mu\medmuskip=2mu plus 1mu minus 2mu\thickmuskip=4mu plus 4mu%
\begin{align*}
P(\mu,\tau)
&=
288\mu^{8}+1048\mu^{7}\tau+1936\mu^{6}\tau^{2}+2688\mu^{5}\tau^{3}
\\
&\phantom{=}
+3012\mu^{4}\tau^{4}+2226\mu^{3}\tau^{5}+941\mu^{2}\tau^{6}+204\mu\tau^{7}+17\tau^{8},
\\
Q(\mu,\tau)
&=
2\mu\tau(144\mu^{7}+524\mu^{6}\tau+968\mu^{5}\tau^{2}+1200\mu^{4}\tau^{3}
\\
&\phantom{=}
+910\mu^{3}\tau^{4}+387\mu^{2}\tau^{5}+84\mu\tau^{6}+7\tau^{7}).
\end{align*}
}

Provided that $\tau=\mu$, the solution takes the form
$$
\lambda
=
P(\mu,\nu)/Q(\mu,\nu),
$$
where
\begin{align*}
P(\mu,\nu)
&=
256\mu^{10}+2112\mu^{9}\nu+8044\mu^{8}\nu^{2}+19355\mu^{7}\nu^{3}
\\
&\phantom{=}
+32167\mu^{6}\nu^{4}+36887\mu^{5}\nu^{5}+28709\mu^{4}\nu^{6}
\\
&\phantom{=}
+14854\mu^{3}\nu^{7}+4912\mu^{2}\nu^{8}+944\mu\nu^{9}+80\nu^{10},
\\
Q(\mu,\nu)
&=
2\mu\nu(\mu+\nu)(192\mu^{8}+1344\mu^{7}\nu+4047\mu^{6}\nu^{2}
\\
&\phantom{=}
+6770\mu^{5}\nu^{3}+6799\mu^{4}\nu^{4}+4216\mu^{3}\nu^{5}+1600\mu^{2}\nu^{6}
\\
&\phantom{=}
+344\mu\nu^{7}+32\nu^{8}).
\end{align*}

\section{Results for Matrices with Arbitrary Constant Entries}

Now we present new results that extend some of the solutions in the previous section to the case of matrices with arbitrary nonnegative constants in place of one or more entries. The results allows one to take account of various possible settings for activity times in actual systems from a random time to constant and zero time.

\subsection{Matrices With One Random Entry}

We examine systems with matrices that have one random entry, whereas the other entries can be arbitrary nonnegative and zero constants. To evaluate the mean cycle time, we apply the same technique as that developed for matrices with zero entries and outlined in the previous section.

The random entry in system matrices is always assumed to have exponential distribution except for one case of a matrix with zero row. In the later case, the particular form of the matrix enables one to obtain a solution that does not rely on exponential distribution assumptions.

\subsubsection{Matrices With Zero Off-Diagonal Entries}

We start with a system that has a matrix
$$
A(k)
=
\left(
\begin{array}{cc}
	\alpha_{k}	& 0 \\
	0	& c
\end{array}
\right),
$$
where $c$ is a nonnegative constant, $\alpha_{k}$ is a r.v. that has exponential distribution with parameter $\mu$. 

The solution procedure that is described in detail in Appendix~\ref{A-MZO} gives the result
$$ 
\lambda
=
c
+
\frac{2e^{-3\mu c}}
{\mu(2-4\mu ce^{-\mu c}+\mu^{2}c^{2}e^{-2\mu c})}.
$$

\begin{figure}[!ht]
\begin{center}
\setlength{\unitlength}{1mm}
\begin{picture}(45,50)

\newsavebox\curvehalfmu
\savebox{\curvehalfmu}(40,40)[b]
{\linethickness{1.5pt}\curve(
2.7,  37.5007,
4.0,  25.0029,
6.0,  16.6789,
8.0,  12.5345,
10.0,  10.0774,
12.0,   8.4817,
14.0,   7.3951,
16.0,   6.6402,
18.0,   6.1141,
20.0,   5.7500,
22.0,   5.5009,
24.0,   5.3325,
26.0,   5.2198,
28.0,   5.1451,
30.0,   5.0958,
32.0,   5.0633,
34.0,   5.0420,
36.0,   5.0279
)}

\newsavebox\curveonemu
\savebox{\curveonemu}(40,40)[b]
{\linethickness{1.5pt}\curve(
2.7,  37.5161,
4.0,  25.0690,
6.0,  16.9633,
8.0,  13.2804,
10.0,  11.5000,
12.0,  10.6649,
14.0,  10.2902,
16.0,  10.1267,
18.0,  10.0559,
20.0,  10.0250,
22.0,  10.0114,
24.0,  10.0053,
26.0,  10.0025,
28.0,  10.0012,
30.0,  10.0006,
32.0,  10.0003,
34.0,  10.0001,
36.0,  10.0001
)}

\newsavebox\curvetwomu
\savebox{\curvetwomu}(40,40)[b]
{\linethickness{1.5pt}\curve(
2.7,  37.8906,
4.0,  26.5608,
6.0,  21.3298,
8.0,  20.2533,
10.0,  20.0500,
12.0,  20.0106,
14.0,  20.0024,
16.0,  20.0006,
18.0,  20.0001,
20.0,  20.0000,
22.0,  20.0000,
24.0,  20.0000,
26.0,  20.0000,
28.0,  20.0000,
30.0,  20.0000,
32.0,  20.0000,
34.0,  20.0000,
36.0,  20.0000
)}

\put(5,5){\vector(1,0){40}}
\put(5,5){\vector(0,1){40}}

\multiput(5,5)(0,10){4}{\line(-1,0){2}}
\multiput(5,5)(10,0){4}{\line(0,-1){2}}

\put(1,0){$0$}
\put(1,14){$1$}
\put(1,24){$2$}
\put(1,34){$3$}
\put(1,44){$\lambda$}

\put(14.5,0){$1$}
\put(24.5,0){$2$}
\put(34.5,0){$3$}
\put(44,0){$\mu$}

\put(5,5){\makebox(40,40){\usebox{\curvehalfmu}}}
\put(5,5){\makebox(40,40){\usebox{\curveonemu}}}
\put(5,5){\makebox(40,40){\usebox{\curvetwomu}}}

\put(29,28){$c=1/2$}
\put(32,17){$c=1$}
\put(32,6.5){$c=2$}

\end{picture}
\caption{Graphs of $\lambda$ as a function of $\mu$.}\label{F-lambda1_mu}
\end{center}
\end{figure}
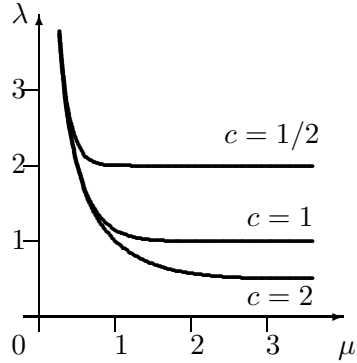

Figures~\ref{F-lambda1_c} and \ref{F-lambda1_mu} show graphs of $\lambda$ as functions of $c$ and $\mu$.
\begin{figure}[!ht]
\begin{center}
\setlength{\unitlength}{1mm}
\begin{picture}(45,50)

\newsavebox\curvehalfc
\savebox{\curvehalfc}(40,40)[b]
{\linethickness{1.5pt}\curve(
0.0,  20.0000,
2.0,  20.0001,
4.0,  20.0023,
6.0,  20.0147,
8.0,  20.0552,
10.0,  20.1549,
12.0,  20.3560,
14.0,  20.7063,
16.0,  21.2487,
18.0,  22.0108,
20.0,  23.0000,
22.0,  24.2038,
24.0,  25.5958,
26.0,  27.1430,
28.0,  28.8125,
30.0,  30.5747,
32.0,  32.4053,
34.0,  34.2855,
36.0,  36.2011
)}

\newsavebox\curveonec
\savebox{\curveonec}(40,40)[b]
{\linethickness{1.5pt}\curve(
0.0,  10.0000,
2.0,  10.0012,
4.0,  10.0276,
6.0,  10.1780,
8.0,  10.6243,
10.0,  11.5000,
12.0,  12.7979,
14.0,  14.4062,
16.0,  16.2026,
18.0,  18.1005,
20.0,  20.0500,
22.0,  22.0251,
24.0,  24.0127,
26.0,  26.0065,
28.0,  28.0033,
30.0,  30.0017,
32.0,  32.0009,
34.0,  34.0005,
36.0,  36.0003
)}

\newsavebox\curvetwoc
\savebox{\curvetwoc}(40,40)[b]
{\linethickness{1.5pt}\curve(
0.0,   5.0000,
2.0,   5.0138,
4.0,   5.3122,
6.0,   6.3989,
8.0,   8.1013,
10.0,  10.0250,
12.0,  12.0063,
14.0,  14.0017,
16.0,  16.0005,
18.0,  18.0001,
20.0,  20.0000,
22.0,  22.0000,
24.0,  24.0000,
26.0,  26.0000,
28.0,  28.0000,
30.0,  30.0000,
32.0,  32.0000,
34.0,  34.0000,
36.0,  36.0000
)}

\put(5,5){\vector(1,0){40}}
\put(5,5){\vector(0,1){40}}

\multiput(5,5)(0,10){4}{\line(-1,0){2}}
\multiput(5,5)(10,0){4}{\line(0,-1){2}}

\put(1,0){$0$}
\put(1,14){$1$}
\put(1,24){$2$}
\put(1,34){$3$}
\put(1,44){$\lambda$}

\put(14.5,0){$1$}
\put(24.5,0){$2$}
\put(34.5,0){$3$}
\put(44,0){$c$}

\put(5,5){\makebox(40,40){\usebox{\curvehalfc}}}
\put(5,5){\makebox(40,40){\usebox{\curveonec}}}
\put(5,5){\makebox(40,40){\usebox{\curvetwoc}}}

\put(7,28){$\mu=1/2$}
\put(7,18){$\mu=1$}
\put(7,7){$\mu=2$}

\end{picture}
\caption{Graphs of $\lambda$ as a function of $c$.}\label{F-lambda1_c}
\end{center}
\end{figure}
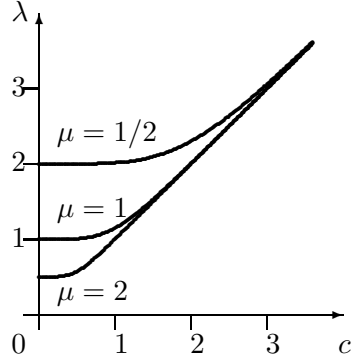

\subsubsection{Matrices With Zero Row or Column}

With similar computations as above, we get the solutions 
$$
A(k)
=
\left(
\begin{array}{cc}
	c	& \beta_{k} \\
	0	& 0
\end{array} 
\right),
\quad
\lambda
=
c
+
\frac{2e^{-2\nu c}}{\nu(2+e^{-2\nu c})},
$$
and
$$
A(k)
=
\left(
\begin{array}{cc}
	c	& 0 \\
	\gamma_{k}	& 0
\end{array} 
\right),
\quad
\lambda
=
c
+
\frac{2e^{-2\sigma c}}{\sigma(2+e^{-2\sigma c})}.
$$

Graphical representation of $\lambda$ with respect to $c$ and $\mu$ is of much the same shape as that of respective graphs depicted in Figures~\ref{F-lambda1_c} and \ref{F-lambda1_mu}.

Now we present a system that admits a solution under more fairly conditions and in a quite general form. Consider a system with a matrix
$$
A(k)
=
\left(
\begin{array}{cc}
	\alpha_{k}	& c \\
	0	& 0
\end{array}
\right),
$$
and assume $c$ to be a nonnegative constant and $\alpha_{k}$ to be a nonnegative r.v. that has a distribution function $F(t)$ with finite mean $a$.

As it is shown in Appendix~\ref{A-MZR}, the mean cycle time for the system is given by
$$
\lambda
=
a
+
\int_{0}^{c}\frac{F(t)F(c-t)(1-F(t))}{1-F(t)F(c-t)}dt.
$$

Specifically, in the case of the exponential distribution function
$$
F(t)
=
\max(0,1-e^{-\mu t}),
$$
we have $a=1/\mu$. Calculation of the integral gives
$$
\lambda
=
\frac{c}{2}+\frac{e^{-\mu c}}{\mu}+\frac{3\arctan\sqrt{4e^{\mu c}-1}-\pi}{\mu\sqrt{4e^{\mu c}-1}}.
$$

\begin{figure}[!ht]
\begin{center}
\setlength{\unitlength}{1mm}
\begin{picture}(55,50)

\savebox{\curvehalfc}(40,40)[b]
{\linethickness{1.5pt}\curve(
0.0,  20.0000,
2.0,  20.0030,
4.0,  20.0220,
6.0,  20.0680,
8.0,  20.1480,
10.0,  20.2663,
12.0,  20.4252,
14.0,  20.6256,
16.0,  20.8674,
18.0,  21.1500,
20.0,  21.4722,
22.0,  21.8326,
24.0,  22.2294,
26.0,  22.6610,
28.0,  23.1256,
30.0,  23.6213,
32.0,  24.1465,
34.0,  24.6994,
36.0,  25.2784,
38.0,  25.8820,
40.0,  26.5086,
42.0,  27.1568,
44.0,  27.8253,
46.0,  28.5129,
48.0,  29.2184
)}

\savebox{\curveonec}(40,40)[b]
{\linethickness{1.5pt}\curve(
0.0,  10.0000,
2.0,  10.0110,
4.0,  10.0740,
6.0,  10.2126,
8.0,  10.4337,
10.0,  10.7361,
12.0,  11.1147,
14.0,  11.5628,
16.0,  12.0732,
18.0,  12.6392,
20.0,  13.2543,
22.0,  13.9127,
24.0,  14.6092,
26.0,  15.3393,
28.0,  16.0990,
30.0,  16.8848,
32.0,  17.6936,
34.0,  18.5228,
36.0,  19.3701,
38.0,  20.2333,
40.0,  21.1107,
42.0,  22.0008,
44.0,  22.9021,
46.0,  23.8134,
48.0,  24.7337
)}

\savebox{\curvetwoc}(40,40)[b]
{\linethickness{1.5pt}\curve(
0.0,   5.0000,
2.0,   5.0370,
4.0,   5.2169,
6.0,   5.5574,
8.0,   6.0366,
10.0,   6.6271,
12.0,   7.3046,
14.0,   8.0495,
16.0,   8.8468,
18.0,   9.6850,
20.0,  10.5554,
22.0,  11.4510,
24.0,  12.3669,
26.0,  13.2987,
28.0,  14.2435,
30.0,  15.1987,
32.0,  16.1622,
34.0,  17.1325,
36.0,  18.1082,
38.0,  19.0885,
40.0,  20.0723,
42.0,  21.0592,
44.0,  22.0484,
46.0,  23.0396,
48.0,  24.0324
)}

\put(5,5){\vector(1,0){50}}
\put(5,5){\vector(0,1){40}}

\multiput(5,5)(0,10){4}{\line(-1,0){2}}
\multiput(5,5)(10,0){5}{\line(0,-1){2}}

\put(1,0){$0$}
\put(1,14){$1$}
\put(1,24){$2$}
\put(1,34){$3$}
\put(1,44){$\lambda$}

\put(14.5,0){$1$}
\put(24.5,0){$2$}
\put(34.5,0){$3$}
\put(44.5,0){$4$}
\put(53,1){$c$}

\put(5,5){\makebox(40,40){\usebox{\curvehalfc}}}
\put(5,5){\makebox(40,40){\usebox{\curveonec}}}
\put(5,5){\makebox(40,40){\usebox{\curvetwoc}}}

\put(7,28){$\mu=1/2$}
\put(10,18){$\mu=1$}
\put(10,7){$\mu=2$}

\end{picture}
\caption{Graphs of $\lambda$ as a function of $c$.}\label{F-lambda2_c}
\end{center}
\end{figure}
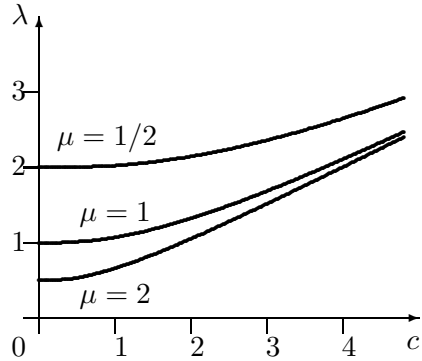

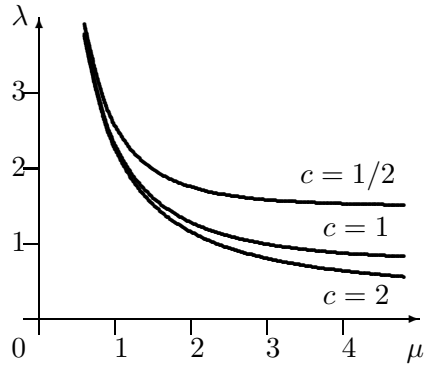
\begin{figure}[!ht]
\begin{center}
\setlength{\unitlength}{1mm}
\begin{picture}(55,45)

\savebox{\curvehalfmu}(40,40)[b]
{\linethickness{1.5pt}\curve(
6.0,  37.5413,
9.0,  25.0850,
12.0,  18.8888,
15.0,  15.1998,
18.0,  12.7658,
21.0,  11.0494,
24.0,   9.7816,
27.0,   8.8125,
30.0,   8.0521,
33.0,   7.4429,
36.0,   6.9467,
39.0,   6.5368,
42.0,   6.1943,
45.0,   5.9053,
48.0,   5.6593
)}

\savebox{\curveonemu}(40,40)[b]
{\linethickness{1.5pt}\curve(
6.0,  37.7776,
9.0,  25.5315,
12.0,  19.5632,
15.0,  16.1042,
18.0,  13.8934,
21.0,  12.3887,
24.0,  11.3187,
27.0,  10.5327,
30.0,   9.9407,
33.0,   9.4859,
36.0,   9.1307,
39.0,   8.8496,
42.0,   8.6245,
45.0,   8.4424,
48.0,   8.2939
)}

\savebox{\curvetwomu}(40,40)[b]
{\linethickness{1.5pt}\curve(
6.0,  39.1265,
9.0,  27.7868,
12.0,  22.6373,
15.0,  19.8814,
18.0,  18.2615,
21.0,  17.2489,
24.0,  16.5878,
27.0,  16.1417,
30.0,  15.8330,
33.0,  15.6151,
36.0,  15.4586,
39.0,  15.3447,
42.0,  15.2609,
45.0,  15.1987,
48.0,  15.1520
)}

\put(5,5){\vector(1,0){50}}
\put(5,5){\vector(0,1){40}}

\multiput(5,5)(0,10){4}{\line(-1,0){2}}
\multiput(5,5)(10,0){5}{\line(0,-1){2}}

\put(1,0){$0$}
\put(1,14){$1$}
\put(1,24){$2$}
\put(1,34){$3$}
\put(1,44){$\lambda$}

\put(14.5,0){$1$}
\put(24.5,0){$2$}
\put(34.5,0){$3$}
\put(44.5,0){$4$}
\put(53,0){$\mu$}

\put(5,5){\makebox(40,40){\usebox{\curvehalfmu}}}
\put(5,5){\makebox(40,40){\usebox{\curveonemu}}}
\put(5,5){\makebox(40,40){\usebox{\curvetwomu}}}

\put(39,23){$c=1/2$}
\put(42,16){$c=1$}
\put(42,7){$c=2$}

\end{picture}
\caption{Graphs of $\lambda$ as a function of $\mu$.}\label{F-lambda2_mu}
\end{center}
\end{figure}

A graphical illustration is given in Figures~\ref{F-lambda2_c} and \ref{F-lambda2_mu}.

\subsubsection{Matrices With Three Constant Entries}

For a system with the symmetric matrix
$$
A(k)
=
\left(
\begin{array}{cc}
	\alpha_{k}	& c \\
	c	& c
\end{array} 
\right),
$$
a solution obtained in \cite{Krivulin2011Calculating} takes the form
$$
\lambda
=
c
+
\frac{2e^{-\mu c}}{\mu(2+e^{-\mu c}-2e^{-2\mu c}+e^{-3\mu c})}.
$$


%
\bibliographystyle{utphys}

\bibliography{Evaluation_of_the_mean_cycle_time_in_stochastic_discrete_event_dynamic_systems}  
%

%
\appendix
\section{Solution Procedures}\label{A-MZO}

\subsection{A Matrix With Zero Off-Diagonal Entries}
We describe evaluation of the mean cycle time for system \eqref{E-zk} with the matrix
$$
A(k)
=
\left(
\begin{array}{cc}
	\alpha_{k}	& 0 \\
	0	& c
\end{array} 
\right),
$$
where $c$ is a nonnegative constant, $\alpha_{k}$ is a r.v. with distribution function
$$
F_{\alpha}(t)
=
\max(0,1-e^{-\mu t}).
$$

Note that equations \eqref{E-xkyk} can now be written as
\begin{align*}
x(k)
&=
\max(x(k-1)+\alpha_{k},y(k-1)),
\\
y(k)
&=
\max(x(k-1),y(k-1)+c).
\end{align*}

With new variables
\begin{align*}
X(k)
&=
x(k)-x(k-1),
\\
Y(k)
&=
y(k)-x(k),
\end{align*}
we arrive at the equations
\begin{align*}
X(k)
&=
\max(\alpha_{k},Y(k-1)),
\\
Y(k)
&=
\max(0,Y(k-1)+c)
-
\max(\alpha_{k},Y(k-1)).
\end{align*}

It is not difficult to see that $X(k)\geq0$, $Y(k)\leq c$, and
$$
x(k)
=
X(1)+\cdots+X(k).
$$

Now the mean cycle time can be represented in the form
\begin{multline*}
\lambda
=
\lim_{k\to\infty}\frac{1}{k}\mathsf{E}\max(x(k),y(k))
\\
=
\lim_{k\to\infty}\frac{1}{k}\mathsf{E}\max(0,Y(k))
+
\lim_{k\to\infty}\frac{1}{k}\sum_{i=1}^{k}\mathsf{E}X(k)
\\
=
\lim_{k\to\infty}\frac{1}{k}\sum_{i=1}^{k}\mathsf{E}X(k).
\end{multline*}

Consider the distribution functions
$$
\Phi_{k}(t)
=
\mathsf{P}\{X(k)<t\},
\qquad
\Psi_{k}(t)
=
\mathsf{P}\{Y(k)<t\},
$$
and note that
$$
\Phi_{k}(t)
=
\mathsf{P}\{\max(\alpha_{k},Y(k-1))<t\}
=
F_{\alpha}(t)\Psi_{k-1}(t).
$$

From the law of total probability it follows that
$$
\Psi_{k}(t)
=
\int_{0}^{\infty}\mathsf{P}\{Y(k)<t|\alpha_{k}=u\}dF_{\alpha}(u).
$$

Evaluation of the conditional probability gives
\begin{multline*}
\mathsf{P}\{Y(k)<t|\alpha_{k}=u\}
\\
=
\mathsf{P}\{\max(0,Y(k-1)+c)-\max(u,Y(k-1))<t\}
\\
=
\begin{cases}
0, & \text{if $u\leq-t$, $t\leq c$}; \\
\Psi_{k-1}(u-c+t), & \text{if $u>-t$, $t\leq c$}; \\
1-\Psi_{k-1}(-t), & \text{if $u\leq-t$, $t>c$}; \\
1, & \text{if $u>-t$, $t>c$}.
\end{cases}
\end{multline*}

After appropriate substitutions, we arrive at the recurrence equation
$$
\Psi_{k}(t)
=
\begin{cases}
\mu e^{\mu t}\int\limits_{0}^{\infty}\Psi_{k-1}(u-c)e^{-\mu u}du, & \text{if $t\leq0$}; \\
\mu\int\limits_{0}^{\infty}\Psi_{k-1}(u-c+t)e^{-\mu u}du, & \text{if $0<t\leq c$}; \\
1, & \text{if $t>c$}.
\end{cases}
$$

Let us investigate the convergence of the sequences $\{\Psi_{k}\}$ and $\{\Phi_{k}\}$. We introduce the notation
$$
a_{k}
=
\mu\int_{-c}^{0}\Psi_{k}(u)e^{-\mu u}du,
\qquad
b_{k}
=
\mu\int_{0}^{c}\Psi_{k}(u)e^{-\mu u}du.
$$

With $C=e^{-\mu c}$, the recurrence equation takes the form 
$$
\Psi_{k}(t)
=
\begin{cases}
C(a_{k-1}+b_{k-1}+C)e^{\mu t}, & \text{if $t\leq0$}; \\
C\left(\frac{c-t}{c}a_{k-1}
+
b_{k-1}
+
C\right)e^{\mu t}, & \text{if $0<t\leq c$}; \\
1, & \text{if $t>c$}.
\end{cases}
$$

It is clear that there is a one-to-one correspondence between the sequence of functions $\{\Psi_{k}\}$ and the sequence of vectors $\{(a_{k},b_{k})\}$. 

The last equation together with the definitions of $a_{k}$ and $b_{k}$ give the recurrences
\begin{align*}
a_{k}
&=
\mu cCa_{k-1}+\mu cCb_{k-1}
+
\mu cC^{2},
\\
b_{k}
&=
\frac{\mu cC}{2}a_{k-1}
+
\mu cCb_{k-1}
+
\mu cC^{2}.
\end{align*}

The iteration matrix of the recurrences has eigenvalues given by
$$
\mu cC\left(1\pm\frac{\sqrt{2}}{2}\right)
=
\mu ce^{-\mu c}\left(1\pm\frac{\sqrt{2}}{2}\right).
$$

Since $0\leq xe^{-x}\leq e^{-1}$ for all $x\geq0$, both eigenvalues have absolute values less than one. Therefore, the recurrences produce iterations that converge to a limiting vector $(a,b)$. The components of the vector are found from the equations
\begin{align*}
(1-\mu cC)a
-
\mu cCb
&=
\mu cC^{2},
\\
-\frac{\mu cC}{2}a
+
(1-\mu cC)b
&=
\mu cC^{2},
\end{align*}
and are given by
$$
a
=
\frac{2\mu cC^{2}}
{2-4\mu cC+\mu^{2}c^{2}C^{2}},
\qquad
b
=
\frac{(2-\mu cC)\mu cC^{2}}
{2-4\mu cC+\mu^{2}c^{2}C^{2}}.
$$

From the convergence of the sequence $\{(a_{k},b_{k})\}$ it follows that the sequence $\{\Psi_{k}\}$ converges as well. The limiting function takes the form
$$
\Psi(t)
=
\begin{cases}
C(a+b+C)e^{\mu t}, & \text{if $t\leq0$}; \\
C(\frac{a}{c}(c-t)+b+C)e^{\mu t}, & \text{if $0<t\leq c$}; \\
1, & \text{if $t>c$}.
\end{cases}
$$

Substitution of the values for $a$ and $b$ yields
$$
\Psi(t)
=
\begin{cases}
\frac{2C^{2}}{2-4\mu cC+\mu^{2}c^{2}C^{2}}e^{\mu t}, & \text{if $t\leq0$}; \\
\frac{2C^{2}}{2-4\mu cC+\mu^{2}c^{2}C^{2}}(1-\mu Ct)e^{\mu t}, & \text{if $0<t\leq c$}; \\
1, & \text{if $t>c$}.
\end{cases}
$$

Furthermore, the sequence $\{\Phi_{k}\}$ also converges to a function given by
\begin{multline*}
\Phi(t)
=
F_{\alpha}(t)\Psi(t)
=
\\
=
\begin{cases}
0, & \text{if $t\leq0$}; \\
\frac{2C^{2}}{2-4\mu cC+\mu^{2}c^{2}C^{2}}(-1+\mu Ct)(1-e^{\mu t}), & \text{if $0<t\leq c$}; \\
1-e^{-\mu t}, & \text{if $t>c$}.
\end{cases}
\end{multline*}

It is easy to see that $\Phi$ presents a distribution function of a r.v. that we denote by $X$. We note that $\mathsf{P}\{X=c\}>0$ whereas $c$ is the only single value the variable $X$ can take with positive probability.

The convergence of the distribution functions $\Phi_{k}\to\Phi$ involves the convergence of means $\mathsf{E}X(k)\to\mathsf{E}X$. This leads to the conclusion that
$$
\lambda
=
\lim_{k\to\infty}\frac{1}{k}\sum_{i=1}^{k}\mathsf{E}X(k)
=
\mathsf{E}X.
$$

To evaluate the last expectation, we first find the probability
\begin{multline*}
p
=
\mathsf{P}\{X=c\}
=
1
-
\int_{0}^{\infty}d\Phi(t)
=
\\
=
1
-
\frac{C(4-6\mu cC+\mu^{2}c^{2}C^{2}-2C+2\mu cC^{2})}{2-4\mu cC+\mu^{2}c^{2}C^{2}}.
\end{multline*}

Now we have
$$
\lambda
=
cp
+
\int_{0}^{\infty}td\Phi(t)
=
c
+
\frac{2C^{3}}
{\mu(2-4\mu c C+\mu^{2}c^{2}C^{2})}.
$$

Finally, substitution $C=e^{-\mu c}$ gives the result
$$
\lambda
=
c
+
\frac{2e^{-3\mu c}}
{\mu(2-4\mu ce^{-\mu c}+\mu^{2}c^{2}e^{-2\mu c})}.
$$

\subsection{A Matrix With Zero Row}\label{A-MZR}

Consider a solution for system \eqref{E-zk} with the matrix
$$
A(k)
=
\left(
\begin{array}{cc}
	\alpha_{k}	& c \\
	0	& 0
\end{array} 
\right),
$$
where $c$ is a nonnegative constant, and $\alpha_{k}$ is a nonnegative r.v. that has a distribution function $F(t)$ with finite mean $a$. First, we represent equations \eqref{E-xkyk} in the form
\begin{align*}
x(k)
&=
\max(x(k-1)+\alpha_{k},y(k-1)+c),
\\
y(k)
&=
\max(x(k-1),y(k-1)).
\end{align*}

By introducing a new state variable
$$
X(k)
=
x(k)-x(k-1),
$$
we reduce \eqref{E-xkyk} to one equation
$$
X(k)
=
\max(\alpha_{k},c-X(k-1)).
$$

Since $x(k)\geq y(k)$ and $x(k)=X(1)+\cdots+X(k)$, we have
$$
\lambda
=
\lim_{k\to\infty}\frac{1}{k}\mathsf{E}\max(x(k),y(k))
=
\lim_{k\to\infty}\frac{1}{k}\sum_{i=1}^{k}\mathsf{E}X(k).
$$

Furthermore, we define the distribution function
$$
\Phi_{k}(t)
=
\mathsf{P}\{X(k)<t\},
$$
and note that
$$
\Phi_{0}(t)
=
\begin{cases}
0, & \text{if $t\leq0$}, \\
1, & \text{if $t>0$},
\end{cases}
\qquad
\Phi_{1}(t)
=
\begin{cases}
0, & \text{if $t\leq c$}, \\
F(t), & \text{if $t>c$}.
\end{cases}
$$

Evaluation of the above probability for all $k\geq1$ gives a recursive equation
$$
\Phi_{k}(t)
=
\begin{cases}
0, & \text{if $t\leq0$}, \\
F(t)(1-\Phi_{k-1}(c-t)), & \text{if $0<t\leq c$}, \\
F(t), & \text{if $t>c$}.
\end{cases}
$$

Let us examine the equation when $0<t\leq c$. After one iteration, we have
$$
\Phi_{k}(t)
=
F(t)(1-F(c-t)(1-\Phi_{k-2}(t))).
$$

Now we define a function
$$
G(t)
=
F(t)F(c-t).
$$

For all even $k=2m$, by iterating the last equation, we arrive at the function
$$
\Phi_{2m}(t)
=
1-(1-F(t))(1+G(t)+\cdots+G^{m-1}(t)).
$$

The sum on the right-hand side is actually an initial part of the geometric series with the common ratio $G(t)<1$.

In the limit as $m$ tends to infinity, for all $t$ under consideration, the subsequence $\Phi_{2m}(t)$ goes to
$$
\Phi(t)
=
1
-
\frac{1-F(t)}{1-G(t)}
=
\frac{F(t)-G(t)}{1-G(t)}.
$$

In the same way, for all $k=2m+1$, we get a subsequence of functions
$$
\Phi_{2m+1}(t)
=
F(t)(1-F(c-t))(1+G(t)+\cdots+G^{m-1}(t)),
$$
and then verify that $\Phi_{2m+1}(t)\to\Phi(t)$ as $m\to\infty$.

From the convergence of the subsequences to a common limit, it follows that when $0<t\leq c$, the entire sequence $\Phi_{k}(t)$ also converges to 
$$
\Phi(t)
=
\frac{F(t)-G(t)}{1-G(t)}
=
\frac{F(t)(1-F(c-t))}{1-F(t)F(c-t)}.
$$

Taking into account that $\Phi_{k}(t)=F(t)$ if $t>c$, we finally conclude that for all $t$, as $k\to\infty$, the sequence $\{\Phi_{k}\}$ converges to the function
$$
\Phi(t)
=
\begin{cases}
0, & \text{if $t\leq0$}, \\
\frac{F(t)(1-F(c-t))}{1-F(t)F(c-t)}, & \text{if $0<t\leq c$}, \\
F(t), & \text{if $t>c$}.
\end{cases}
$$

It is clear that the limiting function $\Phi(t)$ presents a distribution function for a r.v. $X$.

Since $X(k)$ converges to $X$ in distribution, we see that $\mathsf{E}X(k)\to\mathsf{E}X$ as $k\to\infty$, and therefore, $\lambda=\mathsf{E}X$.

To evaluate the expected value of $X$, consider the integral
$$
\int_{0}^{\infty}td\Phi(t)
=
\int_{0}^{c}td\left(\frac{F(t)(1-F(c-t))}{1-F(t)F(c-t)}-F(t)\right)
+
\int_{0}^{\infty}tdF(t).
$$

By applying integration by parts to the first term, we finally get
$$
\lambda
=
a
+
\int_{0}^{c}\frac{F(t)F(c-t)(1-F(t))}{1-F(t)F(c-t)}dt.
$$

\end{document}